\documentclass[12pt]{iopart}
\usepackage{iopams}

\expandafter\let\csname equation*\endcsname\relax   
\expandafter\let\csname endequation*\endcsname\relax
\usepackage{amsmath}

\usepackage{amsfonts}

\usepackage{graphicx} 
\usepackage{xcolor}

\usepackage{subcaption}
\usepackage{tikz}
\usetikzlibrary{calc}
\usepackage{float}
\usepackage{placeins}
\usepackage{multirow}
\usepackage{hyperref}
\hypersetup{
    colorlinks=true,
    linkcolor=blue,
    citecolor=blue,
    urlcolor=blue,
}
\usepackage{cleveref}

\newcommand{\R}{\mathbb{R}}
\begin{document}

\title{Real-time CBCT reconstructions using  Krylov solvers in repeated scanning procedures}

\author{Fred Hastings\textsuperscript{1}\textsuperscript{*}, S M Ragib Shahriar Islam\textsuperscript{2,3,4}\textsuperscript{*}, Malena Sabat\'e Landman\textsuperscript{5}\textsuperscript{*},  Sepideh Hatamikia\textsuperscript{2,3,4}, Carola-Bibiane Sch{\"o}nlieb\textsuperscript{1}, Ander Biguri\textsuperscript{1}}

\address{\textsuperscript{1}Department of Applied Mathematics and Theoretical Physics (DAMTP), University of Cambridge, Cambridge, UK}
\address{\textsuperscript{2}Austrian Center for Medical Innovation and Technology (ACMIT),Wiener Neustadt, Austria}
\address{\textsuperscript{3}Center for Medical Physics and Biomedical Engineering, Medical University of Vienna, Vienna, Austria}
\address{\textsuperscript{4}Danube Private University (DPU), Krems, Austria }
\address{\textsuperscript{5}Department of Mathematics, Emory University, Atlanta, USA}
\address{\textsuperscript{*}These authors contributed equally to this work}

\ead{ander.biguri@gmail.com}
\vspace{10pt}
\begin{indented}
\item[]January 2024
\end{indented}

\begin{abstract}
This work introduces a new efficient iterative solver for the reconstruction of real-time cone-beam computed tomography (CBCT), which is based on the Prior Image Constrained Compressed Sensing (PICCS) regularization and leverages the efficiency of Krylov subspace methods. In particular, we focus on the setting where a sequence of under-sampled CT scans are taken on the same object with only local changes (e.g. changes in a tumour size or the introduction of a surgical tool). This is very common, for example, in image-guided surgery, where the amount of measurements is limited to ensure the safety of the patient. In this case, we can also typically assume that a (good) initial reconstruction for the solution exists, coming from a previously over-sampled scan, so we can use this information to aid the subsequent reconstructions. The effectiveness of this method is demonstrated in both a synthetic scan and using real CT data, where it can be observed that the PICCS framework is very effective for the reduction of artifacts, and that the new method is faster than other common alternatives used in the same setting.

\end{abstract}

%
\vspace{2pc}

\section{Introduction}



Advancements in fast iterative reconstruction techniques have significantly enhanced the utility of Cone Beam Computed Tomography (CBCT) in medical imaging, particularly for interventional procedures that demand rapid image acquisition and reconstruction. These developments are enabling CBCT to transition from primarily diagnostic use to a central role in real-time, image-guided interventions.
Several clinical applications can benefit greatly from such advancements. For example, in image-guided needle-based interventions, such as percutaneous biopsies or robotic-assisted ablations, accurate instrument placement is critical \cite{gulias2020cone, kickuth2015c, hatamikia2023source}. These procedures often require real-time CBCT imaging to guide and verify needle trajectories during the intervention. Since adjustments and repositioning may demand multiple scans within a short time frame, both rapid image reconstruction and low-dose capabilities are essential to maintain procedural efficiency and patient safety. Similarly, spinal surgeries involving pedicle screw placement rely heavily on intraoperative CBCT to verify correct screw positioning \cite{burstrom2021intraoperative, fikuart2025effect}. Misplacement in these contexts can result in irreversible damage, including paralysis. Therefore, high-resolution CBCT imaging is often performed iteratively throughout the procedure. A fast and dose-efficient reconstruction method can significantly reduce surgery time while minimizing the patient's radiation exposure. 
Neuro-interventional procedures, such as cerebral angiography, also frequently involve repetitive CBCT scans. These may be needed to visualize parenchymal structures, monitor hemorrhage progression, or assess the impact of ischemic events \cite{kawauchi2021radioprotection,srinivasan2016cone,irie2008dynact}. As these interventions demand both speed and accuracy, a reconstruction algorithm that delivers high-quality images with fast reconstruction—without compromising dose efficiency—is highly desirable. Furthermore, in image-guided radiation therapy (IGRT), CBCT is routinely used for patient setup and alignment—sometimes daily or multiple times per week where an immediate CBCT scan prior to radiotherapy session is needed for treatment planning. 
Summarizing, two critical challenges emerge across all these applications; first, the need for rapid image acquisition and reconstruction, often in near real-time, to integrate imaging into procedural workflows, second, the need to minimize cumulative radiation dose, especially when multiple scans are required over short time intervals.
Notably, these procedures often involve imaging the same anatomical region repeatedly within a session. This presents a unique opportunity to leverage prior scans as input to the reconstruction algorithm, thereby enhancing image quality even with fewer projections.

Several algorithms have been developed to address the challenges associated with using CBCT in image-guided interventions. One approach involves utilizing reconstruction techniques that inherently reduce radiation exposure compared to the conventional FDK algorithm\cite{feldkamp1984practical}, which is still prevalent in clinical settings. To this end, there has been a recent expansion of iterative solvers that sequentially update the reconstructed image to fit with the measured data for clinical use\cite{scherzer2009variational,liu2014model}. These methods  
have been proven to produce equal quality images under much lower radiation exposure: either by lowering the X-ray dose, or lowering the amount of measurements\cite{baumueller2012low,noel2018evaluation,mileto2019state}. Moreover, these algorithms can easily be adapted to incorporate user-given prior information in the form of regularization. This is often related to assumed properties of the solution, e.g. Tikhonov if the images are smooth or Total Variation (TV) if they are piece-wise constant. More recently, Prior Image Constrained Compressed Sensing (PICCS)\cite{PICCS} was proposed to include structural similarities to previous reconstructions. This type of regularization is very suited for image guided interventions. In fact, assuming a good initial reconstruction, it was shown in \cite{lauzier2012prior} that one can reduce the number of measurements by a very large amount using PICCS regularization in combination with an iterative solver (up to just using 20 projections). However, each iteration of any iterative algorithm takes more time than a single FDK reconstruction, and one generally needs to do tens or often hundreds of iterations to obtain sufficiently a good reconstruction. This has been a roadblock in the clinical adoption of such methods, and is particularly true for the PICCS framework when using the solvers that were proposed in the original and subsequent articles\cite{PICCS,lauzier2012prior,chen2017discriminative,lee2012improved}. To overcome this drawback, we propose using a fast converging class of projection algorithms, namely Krylov subspace methods. These methods converge quadratically (as opposed to the general linear convergence of most standard iterative algorithms) and update the image within the iteration in one step (as opposed to the best performing standard algorithms that require multiple steps per iteration of the whole dataset). In fact, one can reconstruct a clinically sized image within a handful of seconds using these methods, as they require very few iterations to converge\cite{landman2023krylov}. Traditionally, Krylov methods are used to solve linear systems or least-squares problems, which might include Tikhonov regularization, but have also been used in combination with more advanced regularization terms involving $\ell_p$ norms. One very successful approach to do this consists on approximating the $\ell_p$ norms (possibly in combination with a regularization matrix, as in the case of TV regularization), by a sequence of quadratic functionals, giving rise to a chain of (least-squares) minimization problems that need to be solved sequentially \cite{IRN, Daubechies2010IRLS}. For this reason, this technique, known as Iteratively Reweighted Norm (IRN) scheme, involves two nested loops of iterations where the inner loop requires solving a possibly very large-scale least-squares problem. However, this can be done efficiently and with low memory requirements using  Krylov methods based on short-recurrences.  


In this work, we propose a new algorithm aimed at image-guided surgery which combines the benefits of PICCS regularization with the speed of Krylov methods, implemented in an IRN framework. The paper is organized as follows. Section~\ref{sec:methods} concerns the methods of this work. In particular, in Section~\ref{sec:methods1} we describe the mathematical model for CBCT and the PICCS regularization framework. In Section~\ref{sec:methods2} we recall necessary background, and in  Section~\ref{sec:methods3} we describe the new algorithm in detail. Finally, in Section~\ref{sec:results}, we evaluate the performance of the new method 
in simulated data, and subsequently we showcase its use in a lung needle surgery phantom experiment using a Phillips Allura FD20 Xper C-arm CBCT scanner. 

\section{Methods}\label{sec:methods}

This section outlines the methods used in this work. We begin by introducing the problem by formulating the mathematical model for CBCT along with the PICCS regularization framework. Following this, we review the necessary theoretical background for the proposed solution: the IRN scheme; and we finally describe the new algorithm, providing a detailed explanation of its implementation.

\subsection{Prior Image Constrained Compressed Sensing in Krylov methods}\label{sec:methods1} 
Computed Tomography (CT) reconstruction can be modelled as a large scale linear system 
\begin{equation}\label{eq:tomo}
    \mathbf{A}\boldsymbol{x} + \boldsymbol{e} = \boldsymbol{b},
\end{equation}
where \(\boldsymbol{x}\in\R^M\) is the vectorised 3D image with \(M\) voxels, \(\boldsymbol{b}\in\R^N\) is the vector of measurements with added noise, where \(N\) is the number of detector pixels multiplied by the number of projection angles, \(\mathbf{A}\in\R^{N\times M}\) is the sensing matrix and \(\boldsymbol{e}\in\R^N\) is additive noise \cite{kak2001principles}. Even though the noise is not white and Gaussian, this is usually assumed to simplify the computations in experiments that are far from the low photon count limit, see e.g. \cite[Chapter 2.3.2]{Mueller2012Inverse}. When the data is under-sampled or in fast acquisition regimes, recovering the solution $\boldsymbol{x}$ from $\boldsymbol{b}$ in \cref{eq:tomo} can be very sensitive to small perturbations in the measurements. This means that the noise or model errors can easily lead to very noisy reconstructions or artifacts. A good strategy to mitigate this effect is the use of additional information about the reconstruction. In particular, when performing a sequence of CT reconstructions on the same object with only minor differences, we can assume that a previously obtained (and reasonably accurate) reconstruction, denoted as $\boldsymbol{x}_p$, already exists. 
The key idea of this paper is to use information from \(\boldsymbol{x}_p\) to produce accurate reconstructions \(\boldsymbol{x}\) at later times while reducing the X-ray dose. Moreover, to be useful in a clinical setting, this has to be done efficiently in almost real-time.
%

The Prior Image Constrained Compressed Sensing (PICCS) algorithm was proposed in \cite{PICCS} to incorporate information from \(\boldsymbol{x}_p\) in the solution to \cref{eq:tomo} in later times by 
solving the neighbouring variational regularization problem
\begin{equation}
    \min_{\boldsymbol{x}} \left\{\left\lVert\mathbf{A}\boldsymbol{x}-\boldsymbol{b}\right\rVert_2^2 + \alpha^2 \left\lVert\boldsymbol{\Phi}_1(\boldsymbol{x})\right\rVert_1 + \lambda^2 \left\lVert\boldsymbol{\Phi}_2(\boldsymbol{x}-\boldsymbol{x}_p)\right\rVert_1 \right\},
\end{equation}
where \(\boldsymbol{\Phi}_1\) and \(\boldsymbol{\Phi}_2\) are appropriately chosen sparsifying transforms, and where the regularization parameters \(\lambda,\alpha\) are assumed to be known a priori. In particular, they propose choosing $\boldsymbol{\Phi}_i = \boldsymbol{D}$, for $i=1,2$, where $\boldsymbol{D}$ corresponds to a discrete approximation of a differential operator. The resulting regularization terms are also called Total Variation (TV) functionals and the resulting minimization problem becomes:
\begin{equation}\label{PICCS problem}
    \min_{\boldsymbol{x}} \left\{\left\lVert\mathbf{A}\boldsymbol{x}-\boldsymbol{b}\right\rVert_2^2 + \alpha^2 \text{TV}(\boldsymbol{x}) + \lambda^2 \text{TV}(\boldsymbol{x}-\boldsymbol{x}_p)\right\}.
\end{equation}

The first regularization term, $\text{TV}(\boldsymbol{x})$, is known to promote edge reconstructions and therefore reduces the noisy appearance of the recovered solution. Further, practitioners are often most interested in changes between consecutive images; for example observing how a tumour develops over time, or tracking a surgical tool moving in the body. 

Therefore, one might use a $\text{TV}(\boldsymbol{x}-\boldsymbol{x}_p)$ regularization term to promote the salient differences between the two images to be clear. This is because piece-wise constant differences between images have a small TV norm, while small oscillatory differences (that are most likely to correspond to noise) have a large TV norm value. In contrast, penalising $\ell_2$ differences between the current and prior reconstructions, would penalize all differences equally. 

Several variations on the PICCS framework exist in the literature. One such case, which we consider in this paper, is prior-image-registered penalized-likelihood estimation (PIPLE) \cite{Stayman2013PIRPLE}, which involves the following minimization
\begin{equation}\label{PIRPLE_problem}
    \min_{\boldsymbol{x}} \left\{\left\lVert\mathbf{A}\boldsymbol{x}-\boldsymbol{b}\right\rVert_2^2 + \alpha^2 \text{TV}(\boldsymbol{x}) + \lambda^2 \|\boldsymbol{x}-\boldsymbol{x}_p\|_2^2\right\},
\end{equation}
where the regularization term including the prior reconstruction is relaxed to be an $\ell_2$. This is known to have more of an averaging effect between the prior reconstruction and the least-squares solution of the fit-to-data term. Moreover, in practice, a registration step might be necessary to align the prior image $\boldsymbol{x}_p$ with the new reconstructions; this can be seamlessly implemented by adding a suitable deformation to $\boldsymbol{x}_p$ in \eqref{PIRPLE_problem}. This is done, for example, in PIRPLE \cite{Stayman2013PIRPLE}, where a (linear) transformation operator parameterized by a vector of parameters is added to the framework (e.g., modelling 3D rigid body transformation). The framework defined in this paper can be easily adapted to this case, but it is not in the scope of this paper to explore this avenue.

\subsection{The Method}
The strategy used in this paper to solve the minimization problems in \eqref{PICCS problem} and \eqref{PIRPLE_problem} is as follows.  First, the functional in \eqref{PICCS problem} (or \eqref{PIRPLE_problem}) is approximated by a smooth not quadratic functional and, second, this is solved by constructing a sequence of recursive quadratic  approximations (or sub-problems), that  are minimized sequentially. This method, which is part of a wider class of algorithms following a majorization-minimization scheme, requires solving a nested loop of iterations \cite{Lanza2017MM}. The particular approximation strategy that we adopt in this paper, also known as IRN \cite{IRN}, has been used for general TV regularization, and it is known to converge to the minimizer of the smoothed functional. In the following, the mathematical foundations of the presented algorithm are 
reviewed in detail, including the appropriate approximations of the TV, PIPLE, and PICCS functionals, as well as the efficient CGLS solver that we employ for the subproblems in the sequence, which exploits the use of fast and parallelizable GPU applications of the forward model.



\subsubsection{Background on IRN for TV and extension to PIPLE}\label{sec:methods2} 
Let $\mathbf{D}_1 \in \mathbb{R}^{(M^{1/3}-1)\times M^{1/3}}$ be a scaled finite difference approximation of the first derivative operator in one dimension:
$$ \mathbf{D}_1 = 
    \begin{pmatrix}
    1 & -1 & &  \\
     & \ddots &  \ddots&  \\
     &  & 1 & -1      \end{pmatrix}, $$
so that, for 3-dimensional objects, we consider the following discrete approximation of the first derivative operator:
$$  \mathbf{D} =  \begin{pmatrix} \mathbf{D}_x \\ \mathbf{D}_y \\\mathbf{D}_z \end{pmatrix} =
    \begin{pmatrix}
    \mathbf{D}_1 \otimes \mathbf{I} \otimes \mathbf{I}   \\
    \mathbf{I} \otimes \mathbf{D}_1 \otimes \mathbf{I}   \\
    \mathbf{I} \otimes \mathbf{I} \otimes \mathbf{D}_1       \end{pmatrix} ,$$
where $\mathbf{D}_x$,  $\mathbf{D}_y$ and  $\mathbf{D}_z$ correspond to directional derivatives and $\mathbf{I} $ is the identity matrix of dimension $M^{1/3}$.
We are interested in  
the following discrete isotropic total variation (TV) regularization term, which can be defined in 3-dimensional objects as
\begin{equation}\label{TV1}
    \text{TV}(\boldsymbol{x}) =\left\lVert{\sqrt{(\mathbf{D}_x \boldsymbol{x})^2 +(\mathbf{D}_y \boldsymbol{x})^2+(\mathbf{D}_z \boldsymbol{x})^2}}\right\rVert_1,
\end{equation}
where the squaring and square root operations are applied entry-wise. Note that the functional in \eqref{TV1} is not differentiable at any vector with a 0 component, so we consider instead a smoothed version 
defined as
\begin{equation}\label{eq:smooth}
    \text{TV}(\boldsymbol{x}) \approx \left\lVert\mathbf{\widetilde{W}}\left(\mathbf{D}\boldsymbol{x}\right)\mathbf{D}\boldsymbol{x}\right\rVert_2^2,
\end{equation}
with weights defined as
\begin{equation}
    \mathbf{\widetilde{W}}(\mathbf{D}\boldsymbol{x}) = \begin{pmatrix}
        \mathbf{\widetilde{W}_*} & & \\
        & \mathbf{\widetilde{W}_*} & \\
        & & \mathbf{\widetilde{W}_*} 
    \end{pmatrix} \in\mathbb{R}^{3 (M-1)\times M}
\end{equation}
with
\begin{equation}
[\mathbf{\widetilde{W}_*}(\mathbf{D}\boldsymbol{x})]_{ii} = \frac{1}{\sqrt{\sqrt{[\mathbf{D}_x \boldsymbol{x}]_i^2 +[\mathbf{D}_y \boldsymbol{x}]_i^2+[\mathbf{D}_z \boldsymbol{x}]_i^2}+\tau^2}},
\end{equation}
and where $\tau$ is a smoothing parameter fixed ahead of the iterations. 
In this work, this framework is used to find approximate solutions of minimization problems involving TV regularization norms. Assuming that we have an initial guess for the solution $\boldsymbol{x}_0$, then we can define a recursive sequence of reweighted least squares problems of the form
\begin{equation}\label{seq_W}
    \boldsymbol{x}_k = \arg\min_{\boldsymbol{x}}  \left\{\left\lVert\mathbf{A}\boldsymbol{x}-\boldsymbol{b}\right\rVert_2^2 + \alpha^2 \| \mathbf{\widetilde{W}}_{k} \mathbf{D}  \boldsymbol{x}\|_2^2 \right\} \quad \text{ where } \quad \mathbf{{W}}_{k} = \mathbf{\widetilde{W}}(\mathbf{D} \boldsymbol{x}_{k-1}).
\end{equation}
The regularization term in \eqref{seq_W}, ignoring additive constants that do not affect the minimization and multiplicative constants that are absorbed by the regularization parameter, is a quadratic tangent majorant of the smoothed version of the original regularization term, defined in \eqref{eq:smooth}, at $\mathbf{D}\boldsymbol{x}_{k-1}$. By definition, up to constant factors, this means that the regularization term in \eqref{seq_W} is quadratic and is an upper bound for \eqref{eq:smooth}, and that at the point $\mathbf{D}\boldsymbol{x}_{k-1}$, the regularization term in \eqref{seq_W} and \eqref{eq:smooth} take the same value of both the functional and its gradient.
Note that \eqref{seq_W} has a unique solution if $\mathcal{N}(\mathbf{A}) \cap \mathcal{N}(\mathbf{{W}}_{k} \mathbf{D}) \neq \mathbf{0}$, where $\mathcal{N}(\cdot)$ denotes the null-space of a matrix. This is true in CT problems, where the only element in $\mathcal{N}(\mathbf{\widetilde{W}}_{k} \mathbf{D})$ is the constant vector, which is not in the null space of $\mathbf{A}$. This class of methods is known to converge to the solution of the least-squares problem with smoothed TV regularization norm, and it is known as iteratively reweighted least squares (IRLS) or iteratively reweighted norm (IRN) scheme.

Moreover, if one wants to incorporate prior information from a solution $\mathbf{x}_p$, it is also possible to incorporate another term in the minimization problem. For the simplified PIPLE \cite{Stayman2013PIRPLE} with a TV and an $\ell_2$ terms in \eqref{PIRPLE_problem}, this corresponds to:
\begin{equation}\label{seq_W_2}
    \boldsymbol{x}_k = \arg\min_{\boldsymbol{x}}  \left\{\left\lVert\mathbf{A}\boldsymbol{x}-\boldsymbol{b}\right\rVert_2^2 + \alpha^2 \| \mathbf{\widetilde{W}}_{k} \mathbf{D}  \boldsymbol{x}\|_2^2 +\lambda^2 \| \boldsymbol{x}-\boldsymbol{x}_p\|_2^2\right\},
\end{equation}
where $\mathbf{{W}}_{k} = \mathbf{\widetilde{W}}(\mathbf{D} \boldsymbol{x}_{k-1})$.

\subsubsection{Using IRN for PICCS}\label{sec:methods3} 
Analogously to Section \ref{sec:methods2}, we approximate each TV functional in the PICCS problem \eqref{PICCS problem} by a weighted 2-norm (including a mild smoothing):
\begin{equation}\label{weighted2norms}
    \min_{\boldsymbol{x}} \left\{\left\lVert\mathbf{A}\boldsymbol{x}-\boldsymbol{b}\right\rVert_2^2  + \alpha^2 \left\lVert\widetilde{\mathbf{W}}(\mathbf{D}\boldsymbol{x})\mathbf{D}\boldsymbol{x}\right\rVert_2^2+ \lambda^2 \left\lVert\widetilde{\mathbf{W}}(\mathbf{D}\boldsymbol{x}-\mathbf{D}\boldsymbol{x}_p)\mathbf{D}(\boldsymbol{x}-\boldsymbol{x}_p)\right\rVert_2^2\right\},
\end{equation}
with weights that now depend both on the current solution $\boldsymbol{x}$ and on the given reference solution $\boldsymbol{x}_p$. Assuming we have an initial guess for the solution $\boldsymbol{x}_0$, we define the corresponding recursive sequence of reweighted least squares problems as
\begin{equation}\label{PICCS_rec}
    \boldsymbol{x}_k = \min_{\boldsymbol{x}} \left\{\left\lVert\mathbf{A}\boldsymbol{x}-\boldsymbol{b}\right\rVert_2^2  + \alpha^2 \left\lVert\mathbf{W}_k^{(1)}\mathbf{D}\boldsymbol{x}\right\rVert_2^2+ \lambda^2 \left\lVert\mathbf{W}_k^{(2)} \mathbf{D} (\boldsymbol{x}-\boldsymbol{x}_p)\right\rVert_2^2\right\},
\end{equation}
where 
\begin{equation}
    \mathbf{W}_k^{(1)} = \widetilde{\mathbf{W}}(\mathbf{D}\boldsymbol{x}_{k-1}) \quad \text{and} \quad \mathbf{W}_k^{(2)} = \widetilde{\mathbf{W}} (\mathbf{D} \boldsymbol{x}_{k-1}-\mathbf{D} \boldsymbol{x}_p).
\end{equation}

Note that, in practice, problem \eqref{PICCS_rec} is solved by considering the equivalent augmented least-squares problems:
\begin{equation}
    \boldsymbol{x}_k = \min_{\boldsymbol{x}}\left|\left|\begin{bmatrix}
        \mathbf{A}\\
        \alpha \mathbf{W}_k^{(1)} \boldsymbol{D}\\
        \lambda \mathbf{W}_k^{(2)} \boldsymbol{D}
    \end{bmatrix}\boldsymbol{x} - \begin{bmatrix}
        \boldsymbol{b}\\
        \boldsymbol{0} \\
        \lambda\mathbf{W}_k^{(2)} \boldsymbol{D} \boldsymbol{x}_p
\end{bmatrix}\right|\right|_2^2
= \min_{\boldsymbol{x}}\left|\left|\widetilde{\mathbf{A}}\boldsymbol{x}-\Tilde{\boldsymbol{b}}\right|\right|_2^2,
\end{equation}
and approximately solving them using a few inner iterations of a Krylov method (LSQR, or the mathematically equivalent CGLS) on \(\widetilde{\mathbf{A}}\boldsymbol{x}=\tilde{\boldsymbol{b}}\). Note that this is very efficient, since Krylov methods only require matrix-vector products with the system matrix and its transpose, and do not suppose a high memory over cost, since they are based on short recurrences for the solution and residual updates. 

Summarizing, we use the IRN process to solve the TV, PIPLE and PICCS regularization problems so we can leverage the potential of Krylov solvers on large-scale 3D CT problems in an effective way.

\section{Results}\label{sec:results}
In this section, we present two examples representing medical imaging settings where there exists a (reasonably high quality) initial reconstruction, and this is accurate up to local changes, i.e., the introduction of a foreign body in the image, which can represent a tumour or a surgical tool. In particular, the first example corresponds to a simulated head phantom, and associated Digitally Rendered Radiographs (DRR), with an artificial tumour (available in the TIGRE toolbox \cite{biguri2016tigre}). 
The second example corresponds to real scan data, corresponding to an anthropomorphic Thorax phantom with an inserted metallic surgical needle.  Note that the corresponding reconstructions for this example can be subject to strong metal artifacts.

In both examples, we assume that a good initial reconstruction is given, which in practice would most likely correspond to the reconstruction obtained using high-dose measurement data, but that the subsequent scan is heavily undersampled.

In the two examples, we aim to highlight the effect of using different regularization choices, as well as the performance of different optimization algorithms associated to them. Therefore, the reconstructions obtained using the following algorithms are presented, along with other quantitative metrics about their performance:

\begin{itemize}
    \item Direct solver: Feldkamp-Davis-Kress (FDK) algorithm.
    \item Iterative algorithms without explicit regularization: Simultaneous Iterative Reconstruction Technique (SIRT) \cite{SIRT} and Conjugate Gradient Least Squares (CGLS) \cite{CGLS}. Note that most iterative algorithms have implicit regularization properties if equipped with early stopping, see, e.g. \cite{Hansen2010DIP, landman2023krylov}.
    \item Algorithms with only TV regularization on the solution: Iterative Reweighed Norm based Total Variation regularized Conjugate Gradient Least Squares (IRN-TV) \cite{landman2023krylov} and  Adaptive Steepest Descent Projection onto Convex Sets (ASD-POCS-TV) \cite{ASD_POCS}.
    \item Algorithm with TV regularization on the solution and $\ell_2$ norm regularization on the prior term: the proposed Iterative Reweighed Norm based PIPLE algorithm in combination with CGLS (IRN-PIPLE).
    \item Algorithms with PICCS regularization: using ASD-POCS (ASD-POCS-PICCS) \cite{hatamikia2023source} and the proposed algorithm, using an IRN framework in combination with CGLS (IRN-PICCS).
\end{itemize}

The reasoning to compare these algorithms is as follows. First, FDK is the clinical standard. Second, SIRT a standard iterative algorithm, but it displays slow convergence in comparison to the Kylov-based algorithm CGLS. Note that both produce implicitly regularized solutions. Third, ASD-POCS-TV is the most commonly used TV-regularized algorithm in the CT literature, and IRN-TV is a faster Krylov-based algorithm for the same problem. Finally, ASD-POCS-PICCS is the slow converging standard PICCS method, and IRN-PIPLE and IRN-PICCS are the proposed Krylov-based prior image regularized methods.



For all the setups, the corresponding geometry of the scanning devices can be replicated in the TIGRE toolbox when applying the reconstruction algorithms \cite{biguri2016tigre}. 


The PSNR, SSIM and HaarPSI were used for the quantitative analysis of the reconstructed images from all the performed experiments. In particular, HaarPSI is used due to its high correlation with qualitative clinicians' opinion of image quality \cite{Breger_2025}. 

Note that, when we show the reconstructions, only a slice of the 3D volumes is displayed. However, the quantitative error metrics are computed over the whole 3D images.


\subsection{Digital Head Phantom}
\label{sec:HeadImageResults}
In this synthetic example  we consider cases where we want to assess if potential tumours have appeared (or changed in size or shape), using a very low number of projections, and we report reconstruction accuracy and computational time for different choices of regularization and optimization algorithm. In particular, we consider a digital head phantom, of size $128\times128\times128$ voxels, with a (synthetically) added cubic tumour of $8\times8\times8$ voxels. The DRRs projections were computed in TIGRE with resolution $128\times128\times 20$ projections. Note that this is a heavily undersampled dataset: we only consider 20 projections in each imaging step, while a standard clinical scan would require on the order of 600 projections.  Even if this is a simple example, i.e. it is highly unlikely to observe a cubic tumour in practice, we believe this to be a highly illustrative test scenario, in which we can both asses the reconstruction quality of the fine structure in the digital head phantom, as well as the edges of the tumour.

For all the algorithms including prior image regularization, the head phantom image without the tumour was used as prior image. This is a reasonable choice, since we know we can do a high-dose scan once and obtain a very accurate reconstruction. Moreover, since this is a synthetic example, we can use the true solution as the ground truth image in our quantitative analysis. 

Figure \ref{fig:Thead_20p} exhibits a slice of the reconstructed images using all the aforementioned algorithms. Here, (a) is the prior image $\boldsymbol{x}_p$, (f) is the ground truth, and the rests are the reconstructed images (top) with the corresponding reconstructions errors (bottom), all displayed using the same colorbar. Note that the error is just the difference between the reconstruction with respect to the ground truth image, so that high pixel intensity indicates inaccurate reconstructions and a fully black image would indicate a perfect reconstruction. This local error information is more informative than only displaying the error norms, or global error, since one can easily observe there is a high local difference between reconstructions.

\begin{figure*}[h!]
\centering
    \begin{subfigure}[t]{0.9\textwidth}
         \centering
         \includegraphics[width=\textwidth]{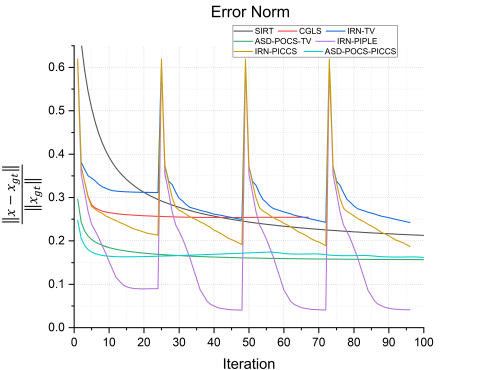}
         \label{fig:exp_1_en_plots_20p}
     \end{subfigure}
\caption{Error norm histories of the digital head phantom reconstructions using 20 projections.}
\label{fig:rn_en_exp1_20p}
\end{figure*}

In terms of global information, Table  \ref{tab:Head_quant_20p} contains the required number of iterations, reconstruction time, and the different image quality metric scores for the reconstructed images. Moreover, Figure~
\ref{fig:rn_en_exp1_20p} shows the relative error norm histories corresponding to the different reconstruction algorithms. 

Since this is a very undersampled dataset, FDK performs qualitatively poorly. Also for SIRT, and IRN-TV the quality drops significantly. The ASD-POCS-TV and ASD-POCS-PICCS are performing well, although a lot of texture information from the image is lost due to the smoothing process. However, the reconstruction time is higher than the other algorithms. The proposed algorithms perform better considering the available texture detail and the reconstruction time. Especially, the image in figure \ref{fig:Head_PICCS01_20p}, IRN-PIPLE, is reconstructing the tumour better by preserving the texture information. Moreover, in figure \ref{fig:Head_PICCS02_20p}, the tumour is more prominently visible due to the added TV regularization in the IRN-PICCS algorithm, however, with compromised original texture information of the image. Note that, for both proposed algorithms, the reconstruction time is in the range of 13 to 76 seconds only.


The error norm histories are displayed in Figure \ref{fig:rn_en_exp1_20p}. 
Note that the methods based on the IRN scheme have inherent restarts (in which the weights for the optimization are updated), which lead to the peaks in the error norms. This is an expected behaviour in these types of methods. However, the proposed IRN-PIPLE method obtains a reconstruction of higher quality in fewer iterations compared to other algorithms. 



\input{imageHead20p.tex}
\FloatBarrier

\begin{table}[]

\begin{tabular}{lccccc}
\hline
\textbf{Algortihm} & \textbf{Iteration} & \textbf{Time} & \textbf{PSNR}    & \textbf{SSIM}   & \textbf{HaarPSI} \\ \hline
FDK            & N/A          & 00:00.2 & 20.6098 & 0.4329 & 0.3827 \\
SIRT           & 100          & 00:30   & 22.2915 & 0.6594 & 0.5152 \\
CGLS           & 20           & 00:04   & 21.5970 & 0.5898 & 0.4448 \\
ASD-POCS       & 100          & 06:57   & 24.5381 & 0.8505 & 0.7091 \\
IRN-TV         & 25           & 00:13   & 20.9140 & 0.6436 & 0.4429 \\
IRN-TV         & 100(4 Outer) & 00:45   & 21.7822 & 0.6895 & 0.4857 \\
IRN-PIPLE      & 25           & 00:18   & 29.4520 & 0.9274 & 0.8524 \\
IRN-PIPLE          & 100(4 Outer)       & 01:09         & \textbf{37.5399} & \textbf{0.9532} & \textbf{0.9394}  \\
IRN-PICCS      & 25           & 00:19   & 22.3996 & 0.7263 & 0.5298 \\
IRN-PICCS      & 100(4 Outer) & 01:16   & 24.1822 & 0.8070 & 0.5710 \\
ASD-POCS-PICCS & 20           & 02:46   & 23.8756 & 0.8012 & 0.6567 \\ \hline
\end{tabular}

\caption{Quantitative results for the digital head phantom example using {20} projections.}
\label{tab:Head_quant_20p}
\end{table}


\FloatBarrier

\FloatBarrier

\subsection{Lung Needle Surgery Scan:}
\label{sec:ThoraxImageResults}

In this example, we reconstruct a Thorax phantom  scanned using the Phillips Allura FD20 Xper C-arm CBCT scanner. This device has source-to-axis and source-to-detector distances of 810 mm and 1195 mm respectively, and an installed detector of size 30$\times$40 cm$^2$ with 0.776 mm pixel pitch. 
Moreover, the data has a resolution of $364\times512\times$\textit{number of projections} and the image is reconstructed with a resolution of $256\times256\times256$ voxels. 

Since the data for this experiment corresponds to real measurements, we consider the ground truth to be a modeled approximation of the true solution. For this, a ground truth of the Thorax phantom is created using 3D Slicer, and the needle data is superimposed artificially in the appropriate location.

The reconstructions corresponding to the real-scan data with different number of projections (180, 50, and 20, respectively), can be observed in Figures \ref{fig:Thorax_180p}, \ref{fig:Thorax_50p}, and \ref{fig:Thorax_20p}. 
In all three figures (a) is the prior image $\boldsymbol{x}_p$ used for prior image regularization, (f) is the ground truth image; and the rest of images correspond to different reconstructions (top) with their corresponding errors (bottom). Analogously to the previous experiment, the difference images (or error) should be completely black in the case of a perfect reconstruction and high intensity indicates more inaccurate reconstructions. 

Note that the real data contains a lot of measurement noise, and due to the presence of the needle, the reconstructions can suffer from strong metal artifacts. Therefore, the reconstructed images using different algorithms, particularly without explicit regularization, tend to show strong errors even with a high number of projections. Accordingly, all the algorithms with explicit regularization perform better in this scenario, and particularly those including PIPLE or PICCS regularization. 
Moreover, the proposed IRN-PIPLE algorithm significantly outperforms other methods, and IRN-PICCS is not far behind. Using only TV regularization struggles to separate artifacts from features, oversmoothing images and thus loosing important image features. Notably, the proposed methods reconstruct high quality images in less than 2 minutes in an implementation that is not optimized to solve this specific geometry.

Some evaluation metrics for this experiment can be seen in Tables \ref{tab:Thorax_quant_180p}, \ref{tab:Thorax_quant_50p}, and \ref{tab:Thorax_quant_20p} (considering 180, 50, and 20 projections respectively) containing the required number of iterations, reconstruction time, and the different image quality metrics scores for the reconstructed images. Note that the image quality measurement metrics are performed over the whole 3D volume with respect to the (approximated) ground truth. Krylov algorithms have been run for 20 iterations, but given the slower convergence of ASD-POCS, 100 iterations are set as maximum. Albeit sometimes these last algorithms will stop earlier due to their stopping criteria, which evaluates the direction of alternating minimization vectors and stops when they are opposite.

The proposed algorithms are the best performing in almost all cases, except for HaarPSI at high number of projections, where the prior image regularization is less meaningful. This can also be qualitatively observed from the reconstructed images.

\input{imageThorax180p.tex}
\FloatBarrier
\begin{table}[!h]
\centering

\begin{tabular}{lccccc}
\hline
\textbf{Algorithm}         & \textbf{Iteration} & \textbf{Time} & \textbf{PSNR} & \textbf{SSIM} & \textbf{HaarPSI} \\ \hline
FDK                        & N/A                & 00:00:01.33    & 33.86         & 0.8202        & 0.5870           \\
SIRT                       & 100                & 01:19          & 40.65         & 0.9464        & 0.6809           \\
CGLS                       & 20                 & 00:15          & 34.05         & 0.8501        & 0.5822           \\
ASD-POCS-TV                   & 100                & 01:15:56       & 36.74         & 0.9463        & 0.6798           \\
IRN-TV                     & 20                 & 01:05          & 37.51         & 0.8935        & 0.6079           \\
IRN-PIPLE                & 20                 & 01:24          & \textbf{43.85 }        & \textbf{0.9473}        & 0.7789           \\
IRN-PICCS                  & 20                 & 01:46          & 41.70         & 0.9181        & 0.6970           \\
ASD-POCS-PICCS             & 100                 & 01:42:22          & 36.74         & 0.8957        & \textbf{0.8434}           \\ \hline
\end{tabular}

\caption{Quantitative results for the thorax phantom example with real measurements involving 180 projections.}
\label{tab:Thorax_quant_180p}
\end{table}
\FloatBarrier

\input{imageThorax50p.tex}
\FloatBarrier

\begin{table}[]
\centering

\begin{tabular}{lccccc}
\hline
\textbf{Algorithm}         & \textbf{Iteration} & \textbf{Time} & \textbf{PSNR} & \textbf{SSIM} & \textbf{HaarPSI} \\ \hline
FDK                        & N/A                & 00:00.62       & 31.02         & 0.6415        & 0.4403           \\
SIRT                       & 100                & 00:55          & 40.68         & 0.9444        & 0.6754           \\
CGLS                       & 20                 & 00:09          & 35.93         & 0.8302        & 0.5748           \\
ASD-POCS-TV                   & 100                & 21:16          & 36.67         & 0.9407        & 0.6626           \\
IRN-TV                     & 20                 & 00:56          & 38.23         & 0.8893        & 0.6152           \\
IRN-PIPLE                  & 20                 & 01:20          & \textbf{44.36}         & \textbf{0.9709}        & \textbf{0.8405}           \\
IRN-PICCS                  & 20                 & 01:43          & 43.76         & 0.9450        & 0.7771           \\
ASD-POCS-PICCS             & 22*                 & 07:33          & 36.92         & 0.9476        & 0.6755           \\ \hline
\end{tabular}

\footnotesize{* set up for 100 iterations but stopped early due to stopping criteria.}
\caption{Quantitative results for the thorax phantom example with real measurements involving 50 projections. }
\label{tab:Thorax_quant_50p}
\end{table}


\FloatBarrier

\input{imageThorax20p.tex}
\FloatBarrier

\begin{table}[]
\centering

\begin{tabular}{lccccc}
\hline
\textbf{Algorithm}         & \textbf{Iteration} & \textbf{Time} & \textbf{PSNR} & \textbf{SSIM} & \textbf{HaarPSI} \\ \hline
FDK                        & N/A                & 00:00.28       & 27.26         & 0.4932        & 0.2971           \\
SIRT                       & 100                & 00:40          & 41.20         & 0.9359        & 0.6795           \\
CGLS                       & 20                 & 00:07          & 37.63         & 0.8301        & 0.6024           \\
ASD-POCS-TV                   & 100                & 08:53          & 37.37         & 0.9385        & 0.6528           \\
IRN-TV                     & 20                 & 00:56          & 40.12         & 0.8935        & 0.6572           \\
IRN-PIPLE                  & 20                 & 01:19          & {43.70}         & \textbf{0.9752}        & 0.8006           \\
IRN-PICCS                  & 20                 & 01:41          & \textbf{43.73}         & 0.9480        & \textbf{0.8158}           \\
ASD-POCS-PICCS             & 19*                 & 02:43          & 38.09         & 0.9462        & 0.6596           \\ \hline
\end{tabular}

\footnotesize{* set up for 100 iterations but stopped early due to stopping criteria.}
\caption{Quantitative results for the thorax phantom example with real measurements involving 20 projections.}
\label{tab:Thorax_quant_20p}
\end{table}

\FloatBarrier

\section{Discussion}

This study discusses different methods to obtain fast and accurate reconstructions in the context of repeated low-dosage CBCT sequence scanning of the same ROI, given a high quality initial scan.
To do so, we focus on both the choice of appropriate regularization to ensure image quality, and the appropriate choice of the algorithm to ensure speed. 

In terms of regularization, we compare implicit regularization given by early stopping of particular iterative methods, total variation (TV) regularization, and methods that combine TV regularization on the image and prior based regularization on the L2 norm (PIPLE) and the TV norm (PICCS). 
First, as noted before in the literature, we confirm that methods based on image prior information (PIPLE and PICCS) very effectively reduce the number of projections required to obtain a meaningful reconstruction. For example, a standard clinical scan for needle surgery requires in the order of 600 projections, or 300 projections in low dose protocols, in each imaging step. This (and other works) suggest instead that using only 3\% of the standard dose can be sufficient, given that a good regularization method is chosen. This reduction is crucial when we consider that, for image-guided needle surgery, several scans are taken during the entire procedure. 


One of the main disadvantages of using advanced regularization is that the most commonly used reconstruction algorithms become significantly more computationally expensive, which can be a bottleneck in their applicability. This is particularly true in image-guided surgery, where very fast reconstructions are needed. In this work, we propose using Krylov methods, a family of iterative methods that show significantly faster convergence than most other commonly used iterative methods.

This is the first time that Krylov methods are used in repeated low-dosage CBCT sequence scanning with image-prior regularization (PIPLE and PICCS). In particular, we combine a technique called iteratively reweighted norm (IRN) approximation, which requires solving a sequence of quadratic subproblems, and we use CGLS as the inner-solver. 
In this work, we observe that these methods significantly outperform other standard algorithms not just in quality but also in computational time for both simulated and real data.

More specifically, we show an example using heavily undersampled data from a digital head phantom which has no or minimum noise component in the projection space. For this synthetic data, we observe that both IRN-PIPLE and IRN-PICCS 
highly perform with only 20 projections, and they have low reconstruction times with respect to the other compared algorithms. Similarly, for an example of a needle surgery dataset with real measurements, IRN-PIPLE and IRN-PICCS provided the best results in the case of having few projections (20 and 50). 

One of the objectives of this work was to provide a PICCS regularized algorithm with a computational speed that would make it clinically viable. Results show that, for clinical data, we can obtain reconstruction times of around 100 seconds, significantly faster than the previously best performing algorithm, ASD-POCS-PICCS. With bespoke implementations in clinical hardware, this could be accelerated further, providing a viable algorithm for on-line image guided interventions. 

It is however important to note that, strictly speaking, ASD-POCS-PICCS and IRN-PICCS (similarly, IRN-TV and ASD-POCS) are solving the exact same mathematical problem, yet this work show their performance differing in practice. This is likely caused by the faster convergence of Krylov methods, reaching a better solution in significantly fewer iterations than ASD-POCS. However, hyperparameter tuning is also a likely factor in this difference. While Krylov methods only have one hyperparameter to set aside from the number of iterations, ASD-POCS-TV has eight, making finding the best performing set a highly complex task\cite{lohvithee2017parameter}. Therefore, even if it is theoretically possible to find a set of parameters and iteration numbers such that both methods perform to the same level of quality, this is not easily done in practice. This is another advantage of Krylov methods: a much simpler set of parameters to set experimentally.

\section{Conclusions}
In this work, we propose two fast and easy to tune algorithms for CBCT imaging when prior images of the same sample/patient are available, using the PIPLE and the PICCS regularization frameworks within a Krylov subspace solver. The resulting algorithms show very good performance, both in image quality and computational time, which makes them very strong algorithmic options for image guided interventions, such as image guided surgery or radiation therapy. 
\FloatBarrier
\section*{Acknowledgements}

 AB acknowledges the support of EPSRC grant EP/W004445/1 and the Accelerate Programme for Scientific Discovery. 
 MSL acknowledges partial support by the U.S. National Science Foundation, Grant DMS-02208294.
 CBS acknowledges support from the Philip Leverhulme Prize, the Royal Society Wolfson Fellowship, the EPSRC advanced career fellowship EP/V029428/1, the EPSRC programme grant EP/V026259/1, and the EPSRC grants EP/S026045/1 and EP/T003553/1, EP/N014588/1, EP/T017961/1, the Wellcome Innovator Awards 215733/Z/19/Z and 221633/Z/20/Z, the European Union Horizon 2020 research and innovation programme under the Marie Skodowska-Curie grant agreement No. 777826 NoMADS, the Cantab Capital Institute for the Mathematics of Information and the Alan Turing Institute. This research was supported by the NIHR Cambridge Biomedical Research Centre (NIHR203312). The views expressed are those of the author(s) and not necessarily those of the NIHR or the Department of Health and Social Care.

\section*{Bibliography} 

\bibliographystyle{unsrt}
\bibliography{main}

\end{document}